\documentclass[12pt,reqno]{amsart}
\newtheorem{theorem}{Theorem}
\newtheorem{proposition}[theorem]{Proposition}
\newtheorem{corollary}[theorem]{Corollary}
\newtheorem{remark}{Remark}

\newfont{\bb}{msbm10 at 12pt}

\def\MIT{\mathrm{MIT}}
\def\lMIT{\lambda^{\MIT}}
\def\APS{\mathrm{APS}}
\def\CHI{\mathrm{CHI}}
\def\mAPS{\mathrm{mAPS}}
\def\O{\Omega}
\def\Si{\Sigma}
\def\N{\mathrm{N}}
\def\gna{\na^{\S}}
\def\D{\mathrm{D}}
\def\L{\mathrm{L}}
\def\R{\mathrm{R}}
\def\H{\mathrm{H}}
\def\P{\mathrm{P}}
\def\SO{\Si\O}
\def\SS{{\bf S}(\pa\O)}
\def\TBO{\mathrm{T}(\pa\O)}
\def\ga{\gamma}
\def\TO{\mathrm{T}\O}
\def\Ga{\Gamma}
\def\A{\mathrm{A}}
\def\im{\mathrm{Im}\,}
\def\BMIT{\mathbb{B}_{\MIT}}
\def\Vpm{\mathrm{V}^{\pm}}
\def\M{\mathrm{M}}
\def\al{\alpha}
\def\S{{\bf S}}

\let\pa\partial     
\let\na\nabla     
     
\DeclareMathAlphabet{\doba}{U}{msb}{m}{n}

\def\End{{\mathop{\rm End}}}

\def\im{{\mathop{\rm Im}}}

\def\Id{{\mathop{\rm Id}}}

\let\<\langle     
\let\>\rangle

\long\def\komment#1{}

\begin{document}


\title[Eigenvalues estimate for the Dirac operator]{Optimal eigenvalues estimate for the Dirac operator on domains with boundary}     

\author{Simon Raulot} 



\date{\today}

\maketitle

\begin{abstract}
We give a lower bound for the eigenvalues of the Dirac operator on a
compact domain of a Riemannian spin manifold under the $\MIT$ bag boundary
condition. The limiting case is characterized by the existence of an
imaginary Killing spinor.
\end{abstract}


\section{Introduction}


Let $\O$ be a compact domain in a $n$-dimensional Riemannian spin manifold
$(\N^n,g)$ whose boundary is denoted by $\pa\O$. In
\cite{hijazi.montiel.roldan:01}, the authors studied four elliptic boundary
conditions for the Dirac operator $\D$ of the domain $\O$. More precisely,
they prove a Friedrich-type inequality \cite{friedrich} which relates the spectrum of the
Dirac operator and the scalar curvature of the domain $\O$. These boundary conditions
are the following: the Atiyah-Patodi-Singer $(\APS)$ 
condition based on the spectral resolution of the boundary Dirac
operator; a modified version of the $\APS$ condition, the $\mAPS$
condition; the boundary condition $\CHI$ associated with a chirality
operator; and a Riemannian version of the $\MIT$ bag boundary condition. In
fact, they show that, if the boundary $\pa\O$ of $\O$ has non-negative mean
curvature, then under the $\APS$, $\CHI$ or $\mAPS$ boundary
conditions, the spectrum of the classical Dirac operator of the domain $\O$
is a sequence of unbounded real numbers $\{\lambda_k:k\in\mathbb{Z}\}$ satisfying
\begin{eqnarray}
\lambda_k^2\geq\frac{n}{4(n-1)}\,\R_0,
\end{eqnarray}

where $\R_0$ is the infimum of the scalar curvature of the domain $\O$. Moreover, equality
holds only for the $\CHI$ and the $\mAPS$ conditions and in these
cases, $\O$ is respectively isometric to a half-sphere or it carries a
non-trivial real Killing spinor and has minimal boundary. In the case of the $\MIT$ boundary condition, 
they show that the spectrum of the Dirac operator on
$\O$ is an unbounded discrete set of complex numbers $\lMIT$ with positive
imaginary part satisfying 
\begin{eqnarray}\label{hmr}
|\lMIT|^2>\frac{n}{4(n-1)}\,\R_0,
\end{eqnarray}

\noindent if the mean curvature of the boundary is non-negative. This result leads to the following question: can one
improve this inequality in order to obtain some boundary geometric invariants on the right hand side of (\ref{hmr})? We show in this paper that such a result can be obtained. More precisely,
we prove the following theorem: 
 
\begin{theorem}\label{main}
Let $\O$ be a compact domain of an $n$-dimensional Riemannian spin manifold
$(\N^n,g)$ whose boundary $\pa\O$ satisfies $\H>0$. Under
the $\MIT$ boundary condition $\BMIT^-$, the spectrum of the classical
Dirac operator $\D$ on $\O$ is an unbounded discrete set of complex numbers
with positive imaginary part. Any eigenvalue $\lMIT$ satisfies
\begin{eqnarray}\label{inemit}
|\lMIT|^2\geq\frac{n}{4(n-1)}\,\R_0+n\,\im(\lMIT)\,\H_0,
\end{eqnarray}
 
where $\H_0$ is the infimum of the mean curvature of the boundary. 
Moreover, equality holds if and only if the associated eigenspinor is an imaginary Killing
spinor on $\O$ and if the boundary $\pa\O$ is a totally umbilical hypersurface with constant mean curvature.
\end{theorem}

The proof of this theorem is based on a modification of the spinorial
Levi-Civita connection which leads to a spinorial Reilly-type
formula. This formula can be seen as a hyperbolic version of the
Reilly inequality used in \cite{hijazi.montiel.roldan:01}.

\noindent The author would like to thank the referee for helpful comments.


\section{Geometric preliminaries}


In this section, we give some standard facts about Riemannian spin
manifolds with boundary. For more details, we refer to \cite{bw:93} or \cite{hijazi.montiel.roldan:01}.

\noindent On a compact domain $\O$ with smooth boundary $\pa\O$ in a $n$-dimensional Riemannian spin manifold $(\N^n,g)$, denote by $\SO$ the complex spinor bundle corresponding to the metric $g$ and by $\nabla$ its Levi-Civita connection acting on $\TO$ as well as its lift to $\SO$. The map $\gamma :\mathbb{C}l\,(\O)\longrightarrow\End(\SO)$ is the Clifford multiplication where $\mathbb{C}l\,(\O)$ is the Clifford bundle over $\O$. The spinor bundle is endowed with a natural Hermitian scalar product, denoted by $\<\,,\,\>$, compatible with $\na$ and $\ga$. The Dirac operator is then the first order elliptic operator acting on sections of $\SO$ locally given by
$$\begin{array}{lcll}
\D: & \Ga(\SO) & \longrightarrow & \Ga(\SO)\\
& \psi & \longmapsto & \sum_{i=1}^{n}\gamma(e_i)\na_{e_i}\psi,
\end{array}$$

where $\{e_1,...,e_n\}$ is a local orthonormal frame of $\TO$. 

\noindent Consider now the boundary $\pa\O$ which is an oriented hypersurface of the domain $\O$ with induced orientation and Riemannian structure. Since the normal bundle of $\pa\O$ is trivial, the boundary itself is a spin manifold. This spin structure on the boundary allows to construct an intrinsic spinor bundle $\Si(\pa\O)$ over $\pa\O$ naturally endowed with a Hermitian metric, a Clifford multiplication $\ga^{\pa\O}$ and a spinorial Levi-Civita connection $\na^{\pa\O}$. Moreover the restriction $\mathbf{S}(\pa\O):=\SO_{|\pa\O}$ to the boundary of the spinor bundle $\SO$ is a Dirac bundle, i.e. there exist on $\mathbf{S}(\pa\O)$ a Hermitian metric denoted by $\<\,,\,\>$ compatible with the Levi-Civita connection $\na^{\mathbf{S}}$ and the Clifford multiplication $\ga^{\mathbf{S}}$. The Clifford multiplication $\ga^{\S}:\mathbb{C}l\,(\pa\O)\longrightarrow\End(\SS)$ is given by $\ga^{\S}(X)\psi=\ga(X)\ga(\nu)\psi$ for all
$X\in\Gamma(\TO)$ and $\psi\in\Gamma\left(\SS\right)$. Similarly we can relate the Levi-Civita connection acting on $\SO$ with that acting on $\mathbf{S}(\pa\O)$ by the spinorial Gauss formula (see \cite{baer:98}):
\begin{eqnarray*}
(\na_X\psi)_{|\pa\O}=\gna_X\psi_{|\pa\O}+\frac{1}{2}\gamma^{\S}(AX)\psi_{|\pa\O},
\end{eqnarray*}

\noindent for all $X\in\Ga\big(\TBO\big)$, $\psi\in\Ga(\SO)$ and where $AX=-\na_X\nu$ is the shape operator of the boundary $\pa\O$ with respect to the inner normal vector field $\nu$. We can then define the boundary Dirac operator acting on $\mathbf{S}(\pa\O)$ which is an elliptic first order differential operator locally given by
\begin{eqnarray}
\D^{\S}=\sum_{j=1}^{n-1}\ga^{\S}(e_j)\na^{\S}_{e_j}.
\end{eqnarray}

\noindent Recall that there is a standard identification
$$\mathbf{S}(\pa\O)\equiv
\left\lbrace
\begin{array}{ll}
\Si(\pa\O) & \quad \text{if n is odd}\\
\Si(\pa\O)\oplus\Si(\pa\O) &  \quad \text{if n is even}
\end{array}\right.$$

\noindent Taking into account the relation between the Hermitian bundle $\mathbf{S}(\pa\O)$ and $\Si(\pa\O)$, one can see that
$$\na^{\mathbf{S}}\equiv
\left\lbrace
\begin{array}{ll}
\na^{\pa\O} & \quad \text{if n is odd}\\
\na^{\pa\O}\oplus\na^{\pa\O} &  \quad \text{if n is even}
\end{array}\right.$$

\noindent and 
$$\ga^{\mathbf{S}}\equiv
\left\lbrace
\begin{array}{ll}
\ga^{\pa\O} & \quad \text{if n is odd}\\
\ga^{\pa\O}\oplus-\ga^{\pa\O} &  \quad \text{if n is even}
\end{array}\right.$$


\section{The $\MIT$ boundary condition}


\noindent First, note that on a closed compact Riemannian spin manifold, the
classical Dirac operator has exactly one self-adjoint $\mathrm{L}^2$ extention,
so it has real discrete spectrum. In the setting of manifolds with
boundary, a defect of self-adjointness appears. It is given by the Green formula
\begin{eqnarray}\label{ipp}
\int_{\O}\<\D\varphi,\psi\>dv(g)-\int_{\O}\<\varphi,\D\psi\>dv(g)=-\int_{\pa\O}\<\ga(\nu)\varphi,\psi\>ds(g),
\end{eqnarray}

\noindent for all $\varphi$, $\psi\in\Ga(\SO)$. Furthermore, in this case, the Dirac operator has a closed range of finite codimension, but an infinite-dimensional kernel, which varies depending on the choice of the Sobolev space. We refer to
\cite{bw:93}, \cite{lo} or \cite{hijazi.montiel.roldan:01} for a careful 
treatment of boundary conditions for elliptic operators. 
 
\noindent The $\MIT$ bag boundary condition has first been introduced by physicists of the
Massachusetts Institute of Technology in a Lorentzian setting (see
\cite{cjjtw}, \cite{cjjt} or \cite{j}). The Riemannian version of this 
condition has been studied in \cite{hijazi.montiel.roldan:01} in order to
get Friedrich estimates and in \cite{hijazi.montiel.zhang:2} because of
its conformal covariance to give a conformal lower bound for the first eigenvalue of the intrinsic Dirac
operator of hypersurfaces bounding a compact domain in a Riemannian spin
manifold. Consider the pointwise endomorphism
\begin{eqnarray*}
i\ga(\nu):\Ga(\SS)\longrightarrow\Ga(\SS)
\end{eqnarray*} 

acting on the restriction to the boundary
$\pa\O$ of the spinor bundle over $\O$ and where $i$ is the fundamental imaginary number. 
This map is an involution, and so the bundle $\SS$ splits into two
eigensubbundles $\Vpm$ associated with the eigenvalues $\pm 1$. We then have
two associated orthogonal projections given by
$$\begin{array}{lccl}
\BMIT^{\pm}: & \L^2(\SS) & \longrightarrow & \L^2(\Vpm)\\ 
 & \varphi & \longmapsto & \frac{1}{2}(\Id\pm i\ga(\nu))\varphi.
\end{array}$$ 

which define local elliptic boundary conditions for the Dirac operator $\D$
on the domain $\O$. So under this boundary condition, the eigenvalue problem 
\begin{equation}\label{evp}
\left\lbrace
\begin{array}{ll}
\D\varphi=\lMIT\varphi & \quad\rm{on}\,\O\\ 
\BMIT^{\pm}\varphi=0 & \quad\rm{along}\,\pa\O
\end{array}
\right.
\end{equation}

has a discrete spectrum with finite dimensional eigenspaces consisting of smooth
spinor fields. 

\begin{remark}\label{rem1}
{\rm  Under the $\MIT$ boundary condition $\BMIT^-$, the spectrum of the
  Dirac operator $\D$ is contained in the upper half complex plane
  $\{z\in\mathbb{C}\,/\,\im(z)>0\}$. Indeed, let $\lMIT$ be an eigenvalue of 
$\D$ under the $\MIT$ boundary condition and $\varphi\in\Ga(\SO)$ the
associated spinor field, then taking $\psi=i\varphi$ in the Formula~(\ref{ipp}) leads to
\begin{eqnarray}\label{immit}
2\,\im(\lMIT)\int_{\O}|\varphi|^2dv(g)=\int_{\pa\O}|\varphi|^2ds(g)
\end{eqnarray}

\noindent Two possibilities can occur: we have either $\im(\lMIT)>0$ or 
$\im(\lMIT)=0$. If $\im(\lambda^{\MIT})=0$, then the spinor field $\varphi$
should vanish along the boundary $\pa\O$ and
by the unique continuation principle (see \cite{bw:93}), it should be
identically zero on the manifold $\O$. This is impossible 
because the spinor $\varphi$ is supposed to be an eigenspinor, so a non
trivial field. The first case is the only possibility,
i.e. $\im(\lMIT)>0$. For the boundary condition $\BMIT^{+}$, we can show that the imaginary part of all eigenvalues of the Dirac operator is negative.}
\end{remark}


\section{The hyperbolic Reilly formula}


In this section, we give a spinorial Reilly formula based on a
modification of the spinorial Levi-Civita connection. Let $\al\in\mathbb{R}$, then
we define the connection $\na^{\al}$ acting on $\SO$ by
\begin{eqnarray}
\na^{\al}_X\varphi:=\na_X\varphi+i\al\ga(X)\varphi,
\end{eqnarray}

for all $\varphi\in\Ga(\SO)$ and $X\in\Ga(\TO)$. We can now derive an
integral version of the Schr\"odinger-Lichnerowicz formula using the 
modified connection $\na^{\al}$. Indeed, we have:
\begin{proposition}
For all spinor fields $\varphi\in\Ga(\SO)$, we have:
\begin{eqnarray}\label{fsl}
\<(\na^{\al})^*\na^{\al}\varphi,\varphi\>_{\L^2}=\<\D^2\varphi,\varphi\>_{\L^2}-\<\frac{\R}{4}\varphi,\varphi\>_{\L^2}+n\al^2||\varphi||^2_{\L^2}-\int_{\pa\O}\<\na^{\al}_{\nu}\varphi,\varphi\>ds(g),
\end{eqnarray}

where $\R$ is the scalar curvature of the domain $\O$.
\end{proposition}

{\it Proof:}
First note that the $\L^2$-formal adjoint of the connection $\na^{\al}$ is,
by definition, given by
\begin{eqnarray*}
\<(\na^{\al})^*\na^{\al}\varphi,\varphi\>_{\L^2}=||\na^{\al}\varphi||^2_{\L^2}=\sum_{j=1}^{n}\int_{\O}\<\na^{\al}_{e_j}\varphi,\na^{\al}_{e_j}\varphi\>dv(g),
\end{eqnarray*}

for all $\varphi\in\Ga(\SO)$ and where $\{e_1,...,e_n\}$ is a local
orthonormal frame of $\TO$. An easy calculation using the compatibility properties of the Hermitian metric with the spinorial connection and the Clifford multiplication gives
\begin{eqnarray*}
\sum_{j=1}^{n}\<\na^{\al}_{e_j}\varphi,\na^{\al}_{e_j}\varphi\>=\sum_{j=1}^{n}\left(e_j\<\na^{\al}_{e_j}\varphi,\varphi\>-\<\na^{-\al}_{e_j}\na^{\al}_{e_j}\varphi,\varphi\>\right),
\end{eqnarray*}

and Stokes theorem leads to
\begin{eqnarray*}
\<(\na^{\al})^*\na^{\al}\varphi,\varphi\>_{\L^2}=\<-\sum_{j=1}^{n}\na^{-\al}_{e_j}\na^{\al}_{e_j}\varphi,\varphi\>_{\L^2}-\int_{\pa\O}\<\na^{\al}_{\nu}\varphi,\varphi\>ds(g).
\end{eqnarray*}

We can now easily compute
\begin{eqnarray*}
\<-\sum_{j=1}^{n}\na^{-\al}_{e_j}\na^{\al}_{e_j}\varphi,\varphi\>_{\L^2} & = & \<-\sum_{j=1}^{n}\na_{e_j}\na_{e_j}\varphi,\varphi\>_{\L^2}+n\al^2||\varphi||^2_{\L^2}\\
& = & \<\na^*\na\varphi,\varphi\>_{\L^2}+n\al^2||\varphi||^2_{\L^2},
\end{eqnarray*}

and then the classical Schr\"odinger-Lichnerowicz formula (see \cite{lawson.michelson:89}) leads to Identity~(\ref{fsl}).
~\hfill$\square$\\

\noindent This formula is a first step to obtain
Inequality~(\ref{inemit}). However, we have now to introduce the Dirac
operator and the twistor operator associated with the connection
$\na^{\al}$. The modified Dirac operator is locally defined by
\begin{eqnarray}
\D^{\al}\varphi=\sum_{j=1}^n\ga(e_j)\na^{\al}_{e_j}\varphi,
\end{eqnarray}

and the associated twistor operator by
\begin{eqnarray}
\P_X^{\al}\varphi=\na^{\al}_X\varphi+\frac{1}{n}\ga(X)\D^{\al}\varphi,
\end{eqnarray}

for all $X\in\Ga(\TO)$ and $\varphi\in\Ga(\SO)$. Note that for $\al=0$, the
operators $\D^{0}$ and $\P^0$ are respectively the classical Dirac operator and
the classical twistor operator which satisfy the relation (see
\cite{bourguignon.hijazi.milhorat.moroianu} or \cite{fried} for example)
\begin{eqnarray*}
|\na\varphi|^2=|\P\varphi|^2+\frac{1}{n}|\D\varphi|^2
\end{eqnarray*} 

We can then check that the modified operators satisfy the same relation,
i.e.
\begin{eqnarray}\label{twist}
|\na^{\al}\varphi|^2=|\P^{\al}\varphi|^2+\frac{1}{n}|\D^{\al}\varphi|^2.
\end{eqnarray} 

Indeed, if $\{e_1,...,e_n\}$ is a local orthonormal frame of $\TO$, we
have
\begin{eqnarray*}
|\P^{\al}\varphi|^2 & = &
\sum_{j=1}^n\<\na^{\al}_{e_j}\varphi+\frac{1}{n}\ga(e_j)\D^{\al}\varphi,
\na^{\al}_{e_j}\varphi+\frac{1}{n}\ga(e_j)\D^{\al}\varphi\>\\
& = & |\na^{\al}\varphi|^2
-\frac{2}{n}|\D^{\al}\varphi|^2+\frac{1}{n}|\D^{\al}\varphi|^2\\
& = & |\na^{\al}\varphi|^2-\frac{1}{n}|\D^{\al}\varphi|^2,
\end{eqnarray*}

and so Identity~(\ref{twist}) follows directly. We are now ready to
establish the hyperbolic version of the spinorial Reilly formula given in
\cite{hijazi.montiel.roldan:01}. This formula can be seen as an 
analogous of the one used in \cite{hijazi.montiel.roldan:02} to give a lower bound of the first
eigenvalue of the intrinsic Dirac operator for hypersurfaces bounding a compact domain of a manifold with negative scalar curvature. More precisely, we
prove:
\begin{proposition}
For all $\varphi\in\Ga(\SO)$, we have:
\begin{eqnarray}\label{frh}
||\P^{\al}\varphi||^2_{\L^2} & = &
\frac{n-1}{n}||\D^{\al}\varphi||_{\L^2}-\<\frac{\R}{4}\varphi,\varphi\>_{\L^2}-n(n-1)\al^{2}||\varphi||^2_{\L^2}\nonumber\\
& & +\int_{\pa\O}\<\D^{\S}\varphi+\frac{n-1}{2}(2\al\,i\ga(\nu)\varphi-\H\varphi),\varphi\>ds(g),
\end{eqnarray}

where $\H$ is the mean curvature of the boundary $\pa\O$ of $\O$.
\end{proposition}

{\it Proof:}
Observe first that the modified Dirac operator $\D^{\al}$ is not formally self-adjoint.
Indeed an easy calculation using~(\ref{ipp}) gives
\begin{eqnarray}\label{ippmod}
\int_{\O}\<\D^{\al}\varphi,\psi\>dv(g)=\int_{\O}\<\varphi,\D^{-\al}\psi\>dv(g)-\int_{\pa\O}\<\ga(\nu)\varphi,\psi\>ds(g),
\end{eqnarray}

for all $\varphi$, $\psi\in\Ga(\SO)$. However, we have:
\begin{eqnarray*}
\D^2\varphi=\D^{-\al}\D^{\al}\varphi-n^2\al^2\varphi,
\end{eqnarray*}

and so substituting in Formula~(\ref{fsl}) gives
\begin{eqnarray*}
\<(\na^{\al})^*\na^{\al}\varphi,\varphi\>_{\L^2}=\<\D^{-\al}\D^{\al}\varphi,\varphi\>_{\L^2}-\<\frac{\R}{4}\varphi,\varphi\>_{\L^2}-n(n-1)\al^2||\varphi||^2_{\L^2}-\int_{\pa\O}\<\na^{\al}_{\nu}\varphi,\varphi\>ds(g).
\end{eqnarray*}

The integration by parts formula~(\ref{ippmod}) leads to
\begin{eqnarray*}
\<(\na^{\al})^*\na^{\al}\varphi,\varphi\>_{\L^2} & = & ||\D^{\al}\varphi||^2_{\L^2}-\<\frac{\R}{4}\varphi,\varphi\>_{\L^2}-n(n-1)\al^2||\varphi||^2_{\L^2}
\\ & & -\int_{\pa\O}\<\ga(\nu)\D^{\al}\varphi+\na^{\al}_{\nu}\varphi,\varphi\>ds(g).
\end{eqnarray*}
 
With the help of Identity~(\ref{twist}), we have
\begin{eqnarray*}
||\P^{\al}\varphi||^2_{\L^2} & = &
\frac{n-1}{n}||\D^{\al}\varphi||_{\L^2}-\<\frac{\R}{4}\varphi,\varphi\>_{\L^2}-n(n-1)\al^{2}||\varphi||^2_{\L^2}\nonumber\\
& & -\int_{\pa\O}\<\ga(\nu)\D^{\al}\varphi+\na^{\al}_{\nu}\varphi,\varphi\>ds(g).
\end{eqnarray*}

However the boundary term can be written
\begin{eqnarray*}
-\ga(\nu)\D^{\al}\varphi-\na^{\al}_{\nu}\varphi =  -\ga(\nu)\D\varphi-\na_{\nu}\varphi+(n-1)\al\,i\ga(\nu)\varphi, 
\end{eqnarray*}

and using the identity 
\begin{eqnarray*}
-\ga(\nu)\D\varphi-\na_{\nu}\varphi=\D^{\S}\varphi-\frac{n-1}{2}\H\varphi,
\end{eqnarray*}

Formula~(\ref{frh}) follows directly.
\hfill$\square$\\

We are now ready to prove Theorem~\ref{main}.


\section{The estimate}


{\it Proof of Theorem~\ref{main}:}
Consider now a compact domain $\O$ of a Riemannian spin manifold such that
the mean curvature $\H$ of the boundary satisfies $\H\geq 2\al$, for
$\al>0$. By ellipticity of the $\MIT$ boundary condition $\BMIT^-$, consider a smooth
spinor field $\varphi\in\Ga(\SO)$ solution of the eigenvalue boundary
problem~(\ref{evp}), i.e. $\varphi$ satisfies
\begin{equation}
\left\lbrace
\begin{array}{ll}\label{mitprob}
\D\varphi=\lMIT\varphi & \rm{on}\,\O\\
\BMIT^-\varphi=0 & \rm{along}\,\pa\O
\end{array}
\right.
\end{equation}

with $\im(\lMIT)>0$ by Remark~\ref{rem1}. We now apply the hyperbolic Reilly
formula~(\ref{frh}) to the spinor field $\varphi$ to get
\begin{eqnarray*}
||\P^{\al}\varphi||^2_{\L^2} & = &
\left(\frac{n-1}{n}|\lMIT-n\al i|^2-n(n-1)\al^{2}\right)||\varphi||_{\L^2}-\<\frac{\R}{4}\varphi,\varphi\>_{\L^2}\\
& & +\int_{\pa\O}\<\D^{\S}\varphi+\frac{n-1}{2}(2\al\,i\ga(\nu)\varphi-\H\varphi),\varphi\>ds(g).
\end{eqnarray*}

Note that since $i\ga(\nu)\varphi=\varphi$ along the boundary, we can compute
$$\<\D^{\S}\varphi,\varphi\>=\<\D^{\S}\varphi,i\ga(\nu)\varphi\>=\<i\ga(\nu)\D^{\S}\varphi,\varphi\>=-\<\D^{\S}\left(i\ga(\nu)\varphi\right),\varphi\>=-\<\D^{\S}\varphi,\varphi\>,$$

and so the preceding formula gives
\begin{align}\label{inmit}
||\P^{\al}\varphi||^2_{\L^2} + \frac{n-1}{2}\int_{\pa\O}(\H-2\al) & |\varphi|^2ds(g)=\\ 
& \frac{n-1}{n}\left(|\lMIT|^2-2n\,\al\,\im(\lMIT)\right)||\varphi||_{\L^2}-\<\frac{\R}{4}\varphi,\varphi\>_{\L^2}\nonumber
\end{align}

The assumption on the mean curvature gives:
\begin{eqnarray*}
|\lMIT|^2-2n\,\alpha\,\im(\lMIT)\geq\frac{n}{4(n-1)}\,\R_0.
\end{eqnarray*}

For $\al_{0}=\frac{1}{2}\,\H_0$, where $\H_0=\inf_{\pa\O}(\H)$, we get Inequality (\ref{inmit}). 
Suppose now that equality is achieved, thus
\begin{eqnarray*}
||\P^{\al_0}\varphi||^2_{\L^2}=0\quad\textrm{and}\quad\frac{n-1}{2}\int_{\pa\O}(\H-2\al_0)|\varphi|^2ds(g)=0.
\end{eqnarray*}
 
Moreover the spinor field $\varphi$ is a solution of (\ref{mitprob}), so it satisfies the Killing equation
\begin{eqnarray*}
\na_{X}\varphi=-\frac{\lMIT}{n}\ga(X)\varphi,\qquad\text{for all}\;\;X\in\Ga(\TO).
\end{eqnarray*} 

Since such a spinor field has no zeroes (see \cite{fried}), the mean
curvature of the boundary is constant with $\H=2\al_0$. Furthermore, it is a well-known
result \cite{baum.friedrich.grunewald.kath:90} that, in this case, the
eigenvalue $\lMIT$ has to be either real or purely imaginary. Here we have $\im(\lMIT)>0$, then $\lMIT\in
i\mathbb{R}^+_*$. The domain $\O$ is in particular an Einstein
manifold. We now show that the boundary has to be totally umbilical. Indeed, note that we have for all $X\in\Ga(\TBO)$:
\begin{eqnarray*}
\na_X(i\ga(\nu)\varphi) & = & i\ga(\na_X\nu)\varphi+i\ga(\nu)\na_{X}\varphi\\
 & = & i\ga(\na_X\nu)\varphi+\al_0\ga(\nu)\ga(X)\varphi\\
& = & i\ga(\na_X\nu)\varphi-\al_0\ga(X)\ga(\nu)\varphi\\
& = & i\ga(\na_X\nu)\varphi+i\al_0\ga(X)\varphi.
\end{eqnarray*} 
\noindent However along the boundary we have $i\ga(\nu)\varphi=\varphi$, so we obtain
\begin{eqnarray*}
\ga(\na_X\nu)\varphi=-2\al_0\ga(X)\varphi.
\end{eqnarray*} 

\noindent Since the spinor field $\varphi$ has no zeros, we have $\A(X)=-\na_X\nu=2\al X$ and the boundary is 
totally umbilical.We can again show that in the equality case, we have $\im(\lMIT)=n\al_0$. In fact, just note
that the boundary term can be rewritten as
\begin{eqnarray*}
\int_{\pa\O}\<\D^{\S}\varphi-\frac{n-1}{2}\H\varphi+(n-1)\al_0\varphi,\varphi\>ds(g)=-\int_{\pa\O}\<\na_{\nu}\varphi+\ga(\nu)\D\varphi-(n-1)\al_0\varphi,\varphi\>ds(g).
\end{eqnarray*}

This term is zero since we have equality in (\ref{inmit}). Now using that the spinor field $\varphi$ is an imaginary Killing spinor satisfying (\ref{evp}) gives
\begin{eqnarray*}
\na_{\nu}\varphi+\ga(\nu)\D\varphi=\frac{n-1}{n}\im(\lMIT)\varphi.
\end{eqnarray*}

Substituting in the preceding identity gives
\begin{eqnarray*}
(n-1)\int_{\pa\O}(\al_0-\frac{\im(\lMIT)}{n})|\varphi|^2ds(g)=0,
\end{eqnarray*}

and since $\varphi$ has no zeroes, $\im(\lMIT)=n\al_0=\frac{n\H_0}{2}$.
\hfill$\square$\\

\begin{remark}\label{rem2}

\noindent 

\begin{enumerate}

\item {\rm The orthogonal projection $\BMIT^+$ defines a local elliptic boundary
    condition for the Dirac operator $\D$ of $\O$. We can easily check
    that in this case, the imaginary part of an eigenvalue $\lMIT$ of $\D$
    satisfies $\im(\lMIT)<0$. Inequality~(\ref{inemit}) is then given by}
\begin{eqnarray*}
|\lMIT|^2\geq\frac{n}{4(n-1)}\,\R_0-n\,\im(\lMIT)\,\H_0.
\end{eqnarray*} 

\item {\rm For $\H_0=0$, we obtain Inequality~(\ref{hmr}). In fact, if we suppose that equality is achieved, Theorem~\ref{main} implies $\im(\lMIT)=\frac{n\H_0}{2}=0$ which is impossible by Remark~\ref{rem1}.}

\item  {\rm Note that the Riemannian spin manifolds with an imaginary Killing spinor with Killing number $i\al$ have
  been classified by H. Baum in \cite{h1} and \cite{h2}. Such manifolds are
  called pseudo-hyperbolic and they are given by
\begin{eqnarray*}
(\mathbb{R}\times_{\exp}\M_0,g)=(\mathbb{R}\times\M_0,dt^2\oplus e^{-4\al t}g_{\M_0}),
\end{eqnarray*}

\noindent where $(\M_0,g_{\M_0})$ is a complete Riemannian spin manifold carrying a
non-trivial parallel spinor. After suitable rescaling of the metric, we can
assume that the Killing number is either $i/2$ or $-i/2$, i.e. we have
\begin{eqnarray*}
\na_{X}\phi=\pm\frac{i}{2}\ga(X)\phi.
\end{eqnarray*}

\noindent Moreover, constant mean curvature 
hypersurfaces in pseudo-hyperbolic manifolds are classified by the Hyperbolic Alexandrov
Theorem proved in \cite{montiel} (see also \cite{hijazi.montiel.roldan:02} 
for a proof using spinors). Indeed, such a hypersurface is either a round
geodesic hypersphere (and, in this case, $\M_0$ is flat and $\H>1$) or a
slice $\{s\}\times\M_0$ (and, in this case, $\M_0$ is compact and
$\H=1$).}
\end{enumerate}
\end{remark}

We can then prove the following corollary:
\begin{corollary}
If the boundary of the compact domain $\O$ is connected, there is no manifold satisfying the equality case in Inequality~(\ref{inemit}).
\end{corollary} 

{\it Proof:} If $\O$ is a compact domain whith connected boundary achieving equality in~(\ref{inemit}), then there exists an imaginary Killing spinor on $\O$ and the boundary $\pa\O$ is a totally umbilical constant mean curvature hypersurface with $\H=2\al$. However, using Remark~(\ref{rem2})$.3$, $\O$ is a domain in a pseudo-hyperbolic space whose connected boundary is a slice $\{s\}\times\M_0$ and then $\O$ is non-compact.
\hfill $\square$\\

\begin{remark}
\rm{ With a slight modification of the boundary condition, we give a domain $\O$ whose boundary has two connected components carrying an imaginary Killing spinor field $\varphi\in\Ga(\SO)$ which satisfy
\begin{eqnarray}\label{twobound}
i\ga(\nu_1)\varphi_{|\pa\O_1}=\varphi_{|\pa\O_1}\qquad\text{and}\qquad i\ga(\nu_2)\varphi_{|\pa\O_2}=-\varphi_{|\pa\O_2},
\end{eqnarray}

\noindent where $\nu_1$ (resp. $\nu_2$) is an inner unit vector field normal to $\pa\O_1$ (resp. $\pa\O_2$). First recall that one distinguishes two types of imaginary Killing spinors (see \cite{h1} and \cite{h2}). Indeed, if $\varphi\in\Gamma(\SO)$ is an imaginary Killing spinor, denote by $f$ its length function, then the function
\begin{eqnarray*}
q_{\varphi}(x):=f(x)^2-\frac{1}{4\al^2}||\na f||^2
\end{eqnarray*}

\noindent satisfies $q_{\varphi}$ is constant and $q_{\varphi}\geq 0$. If $q_{\varphi}=0$, $\varphi$ is a Killing spinor of type I whereas if $q_{\varphi}>0$, $\varphi$ is a Killing spinor of type II. If $(\N^n,g)$ is a complete connected Riemannian spin manifold with an imaginary Killing spinor of type II associated with the Killing number $i\al$, then $(\N^n,g)$ is isometric to the hyperbolic space $\mathbb{H}^n_{-4\al^2}$. If $(\N^n,g)$ admits an imaginary Killing spinor of type I, then $(\N^n,g)$ is isometric to the warped product $(\mathbb{R}\times\M_0,dt^2\oplus e^{-4\al t}g_{\M_0})$, where $\M_0$ is a complete Riemannian spin manifold with a non-trivial parallel spinor field. Moreover, $q_{\varphi}=0$ if and only if there exists a unit vector field $\xi$ on $\N$ such that $\ga(\xi)\varphi=i\varphi$. In fact, we can easily prove that the vector field $\xi$ is the normal field of $\{t\}\times\M_0$ for all $t\in\mathbb{R}$. So consider the domain given by the warped product 
$\O:=([a,b]\times\M_0,dt^2\oplus e^{-4\al t}g_{\M_0}),$ where $\M_0$ is a compact spin manifold carrying a non-trivial parallel spinor field and with $-\infty<a<b<+\infty$. The domain $\O$ carries an imaginary Killing spinor $\varphi$ of type I, so there exists $\xi$ normal to $\{t\}\times\M_0$ for all $t\in [a,b]$ such that $\ga(\xi)\varphi=i\varphi$. The boundary of $\O$ has two connected components which are slices $\{a\}\times\M_0$ and $\{b\}\times\M_0$ of $\O$ and with mean curvature $\H_a=\H_b=2\al$, where $\H_t$ is the mean curvature of a slice $\{t\}\times\M_0$. The spinor field $\varphi$ clearly satisfies the boundary conditions (\ref{twobound}).  
}
\end{remark}


\bibliographystyle{amsalpha}     
\bibliography{bibnote}


\vspace{1cm}     
Author address:     
\nopagebreak     
\vspace{5mm}\\     
\parskip0ex     
\vtop{     
\hsize=6cm\noindent     
\obeylines     

}     
\vtop{     
\hsize=8cm\noindent     
\obeylines     
Simon Raulot,
Institut \'Elie Cartan BP 239     
Universit\'e de Nancy 1     
54506 Vand\oe uvre-l\`es -Nancy Cedex     
France     
}     
     
\vspace{0.5cm}     
     
E-Mail:     
{\tt raulot@iecn.u-nancy.fr  }

\end{document}